\documentclass[10pt]{article}
\usepackage[utf8]{inputenc}
\usepackage{authblk}
\usepackage[style=numeric-comp,sorting=none]{biblatex}

\usepackage[a4paper, margin=3.5cm]{geometry}
\usepackage{listings}
\usepackage{graphicx}
\usepackage{verbatim}
\usepackage{stackengine}
\usepackage{subcaption}
\usepackage{siunitx}
\usepackage{xspace}
\usepackage{amsmath}
\usepackage{tablefootnote}
\usepackage{mathtools}

\begin{document}

\title{\textbf{On the sum of ordered spacings}}
\author{Lolian Shtembari}
\author{Allen Caldwell}
\affil{Max Planck Institute for Physics, Munich}
\date{July 2020}

\maketitle

\begin{abstract}

\noindent We provide the analytic forms of the distributions for the sum of ordered spacings. We do this both for the case where the boundaries are included in the calculation of the spacings and the case where they are excluded.  Both the probability densities as well as their cumulatives are provided.  These results will have useful applications in the physical sciences and possibly elsewhere.

\end{abstract}

\section{Introduction}
\label{sec:introduction}

The use of spacings between ordered real-valued numbers is very useful in many areas of science. In particular, either unnaturally small or large spacings can be a signal of an interesting effect. As particle physicists, we are interested in the appearance of the unexpected clustering of values, indicating the presence of a new process, or large gaps between the ordered values, allowing us to set upper limits on the normalization of a distribution. 


Order statistics have been studied in great depth in the statistics community, but the work is poorly known in the physics community.  This has led to the rediscovery of results long-known to the statistics community.  An example of the use of spacings between values in the particle physics community is presented by Yellin \cite{Yellin:2002}, where the author proposes a method to set a limit on the interaction rate of putative dark matter particles using the size of gaps in the observed energy spectrum of recorded interactions. In this context, a large gap in the energy spectrum implies an upper limit on the interaction strength.  The `Yellin method' has been used by many groups in reporting their results (see for example~\cite{ref:CRESST}).

As is normally done when using order statistics, a known distribution with some quantity described as a function of some variable is converted to the uniform distribution $\mathcal{U} \in [0,1]$ via the cumulative of the distribution of interest.  The powerful results that can be derived for order statistics on the unit interval can then be applied to the task at hand.  The maximal spacing distribution given by Yellin was previously known for at least 75 years~ \cite{Pincus:1984} .

 In the following we will introduce some of the main concepts at the root of order statistics.  We will then introduce new statistics formed from the sum of ordered spacings.  We believe these statistics can be of use in scientific applications, and, as far as we know, have not been presented to date. For a more comprehensive and detailed review of the results that have been achieved in the field of order statistics, we recommend the following books \cite{DavidNagaraja:2003} \cite{Nagaraja:2008}.\\
\\
\section{Notation and Known Results}
Let $\{X_1, X_2,..., X_N\}$ be a sequence of iid random variables with a uniform distribution in the range $[0,1]$. Given this sample of variables, its \textit{k}-th order statistic is defined as the \textit{k}-th smallest value of the sample. This means that given the sequence $\{X_1, X_2, ..., X_N\}$ we obtain the order statistic $\{X_{(1)}, X_{(2)}, ..., X_{(N)}\}$.

\begin{equation}
    X_{(k)} \coloneqq \text{ the \textit{k}-th smallest of $\{X_1, X_2, ..., X_N\}$ for $k = 1,2,...,N$}
\end{equation}
\\

\noindent The distribution of these ordered samples is simply a Beta distribution:

\begin{equation} \label{os_dist}
    X_{(k)} \sim Beta(k, N + 1 - k)
\end{equation}
\\
and the joint distribution of two ordered values is:

\begin{equation} \label{gos_dist}
    X_{(k)} - X_{(j)} \sim Beta(k - j, N + 1 - (k - j))
\end{equation}
\\
Given $\{X_{(1)}, X_{(2)}, ..., X_{(N)}\}$ we are interested in the the spacings between the ordered values. 

To begin with, we will consider an extended set of ordered values, namely the boundaries of the range of $X_{(i)}$ themselves: 0 and 1. We will define $X_{(0)} = 0$ and $X_{(N+1)} = 1$. We define a spacing $G_i$ as the distance between neighboring values:

\begin{equation} 
    G_i = X_{(i)} - X_{(i-1)} \text{   for $i = 1,2,...,N+1$}
\end{equation}
\\
The distribution of any of these is given by Eq.~(\ref{os_dist})):

\begin{equation} \label{pdf_g1}
    p(G_i = x) = N(1 - x)^{N - 1}
\end{equation}
\begin{equation} \label{cdf_g1}
    p(G_i \leq x) = 1 - (1 - x)^N
\end{equation}
\\
Just as we ordered the set of initial values $\{X_1, X_2$, ... ,$X_N\}$ in order to obtain $\{X_{(0)}, X_{(1)}$, ... ,$X_{(N+1)}\}$, we can order the set $\{G_1, G_2$, ... ,$G_{N+1}\}$ and obtain $\{G_{(1)}, G_{(2)}$, ... ,$G_{(N+1)}\}$ which could be interpreted as an order statistics of the spacings.
\\

\noindent The distribution of the smallest spacing $G_{(1)}$, due to R.~Fisher~\cite{Fisher:1929}, is known:

\begin{equation}\label{k-th-spacing}
    p(G_{(1)} = x | N, 1) = N (N+1) \left(1 - (N+1)x  \right)^{N-1}
\end{equation}

\noindent as is the general distribution of $G_{(k)}$ given $N$ uniform samples in an interval of length $1$ \cite{Feller:1966}:

\begin{equation}\label{k-th-spacing}
    p(G_{(k)} = x | N, 1) = N (N+1) \binom{N}{k-1} \sum_{i=1}^k (-1)^{k-i}  \binom{k-1}{i-1} \left[ 1 - (N+2-i)x \right]^{N-1} \cdot H \left(x, 0, \frac{1}{N+2-i} \right)
\end{equation}
\\
where $H(x, a, b) = 1$ if $a \leq x \leq b$ and 0 otherwise. Limits of this distributions have also been evaluated for different combinations of $k$ and $N$ \cite{Holst:1980, Bairamov:2010}.
\\

\noindent Given this new set of ordered uniform spacings, an interesting quantity is the sum of the first $k$ smallest or largest $G_{(i)}$.  We will derive the probability distributions for these quantities in the following.
\\

\noindent We denote the sum of the first $k$ minima as:

\begin{equation} \label{def_sum_minima}
    s_k = \sum_{i = 1}^{k} G_{(i)}
\end{equation}

\noindent and the sum of the first $k$ maxima as

\begin{equation} \label{def_sum_minima}
    S_k = \sum_{i = 1}^{k} G_{(N+2-i)}
\end{equation}
\section{Sum of smallest spacings}

\subsection{Sum of minima: $k=1$}

If $k = 1$ then $s_k = G_{(1)}$, so the distribution of $s_1$ is the same as that of the smallest spacing:

\begin{equation} \label{pdf_min1}
    p(s_1 = s) = (N+1)N \left[ 1 - (N+1)s \right]^{N - 1}
\end{equation}
\begin{equation} \label{cdf_min1}
    p(s_1 \leq s) = 1 - \left[ 1 - (N+1)s \right]^N
\end{equation}

\subsection{Sum of minima: $k=2$}

In order to get the distribution of $s_2$, it is useful to consider the joint distribution of $(G_{(1)}, s_2)$. Using the Law of Total Probability:

\begin{equation} \label{min1_sum2_bayes_decomposition}
    p(G_{(1)} = x, s_2 = s | N, 1) = p(G_{(1)} = x | N, 1) \cdot p(s_2 = s | N, 1, G_{(1)} = x)
\end{equation}
\\
In order to derive an expression for $p(s_2|G_{(1)})$ we can consider that once we have chosen the length on the smallest spacing, by definition all the other spacings need to be longer or equal to this minimum length. We can then proceed to subtract $G_{(1)}$ from the length of all the other spacings:

\begin{equation} \label{rescale_spacings}
    G_{(i)} - G_{(1)} = G_{(i-1)}^* \text{ , for } i = 2, ..., N+1 
\end{equation}
\\
This operation leaves us with a reduced set of spacings (since subtracting $G_{(1)}$ from itself results in 0, we simply discard this element) where the reduced spacings are still sorted in increasing order:
\begin{equation} \label{rescaled_spacings_set}
    \{ G_{(1)}, ..., G_{(N+1)} \} \rightarrow \{ G_{(1)}^*, ..., G_{(N)}^* \}
\end{equation}
\\
and they sum up to:
\begin{equation} \label{rescaled_spacings_sum}
    \sum_{i = 1}^N G_{(i)}^* = 1 - (N+1) G_{(1)}
\end{equation}
\\
The set $\{G_{(1)}^*$, ... ,$G_{(N)}^*\}$ can be interpreted as ordered uniform spacings determined by sampling $N-1$ values in an interval of length $1 - (N+1) G_{(1)}$.
Given this rearrangement, we can express the sum of $k$ minima using this new set of  spacings:
\begin{equation} \label{rescale_sum_minima}
    s_{k-1}^* = \sum_{i=1}^{k-1} G_{(i)}^* =  \sum_{i=1}^{k} (G_{(i)} - G_{(1)}) = s_k - k \cdot G_{(1)}
\end{equation}
\\
This allows us to rewrite the conditional distribution of $s_2$ as:

\begin{align} \label{pdf_sum2_given_min1}
    p(s_2 = s | N, 1, G_{(1)} = x) &= p(s_1^* = s - 2 x | N-1, 1 - (N + 1)x) \nonumber \\
    &= \left( \frac{1}{1 - (N+1)x} \right) N (N-1) \left[ 1 - \frac{N(s - 2x)}{1 - (N+1)x} \right]^{N-2}
\end{align}
\\
Putting Eq.~(\ref{pdf_min1}, \ref{min1_sum2_bayes_decomposition}, \ref{pdf_sum2_given_min1}) together we obtain:

\begin{align} \label{pdf_min1_sum2}
    p(G_{(1)} = x, s_2 = s) &= (N+1) N^2 (N-1) \left[ 1 - (N+1)x \right]^{N-2} \left[ \frac{1 + (N-1)x - N s}{1 - (N+1) x} \right]^{N-2} \nonumber \\
    &= (N+1) N^2 (N-1) \left[ 1 - (N - 1)x - N s \right]^{N-2}
\end{align}
\\
The support of $s_2$ is $\left[ 0, \frac{2}{N + 1} \right]$ and the support of $s_1^*$ with $N-1$ samples is $\left[ 0, \frac{1}{N} \right]$, thus the joint distribution is bound withing a triangle as showed in Fig.~{\ref{min1_sum2}}. Marginalizing over $G_{(1)}$ we get the distribution of $s_2$:

\begin{figure} [h!]
    \centering
    \includegraphics[width=0.5\textwidth]{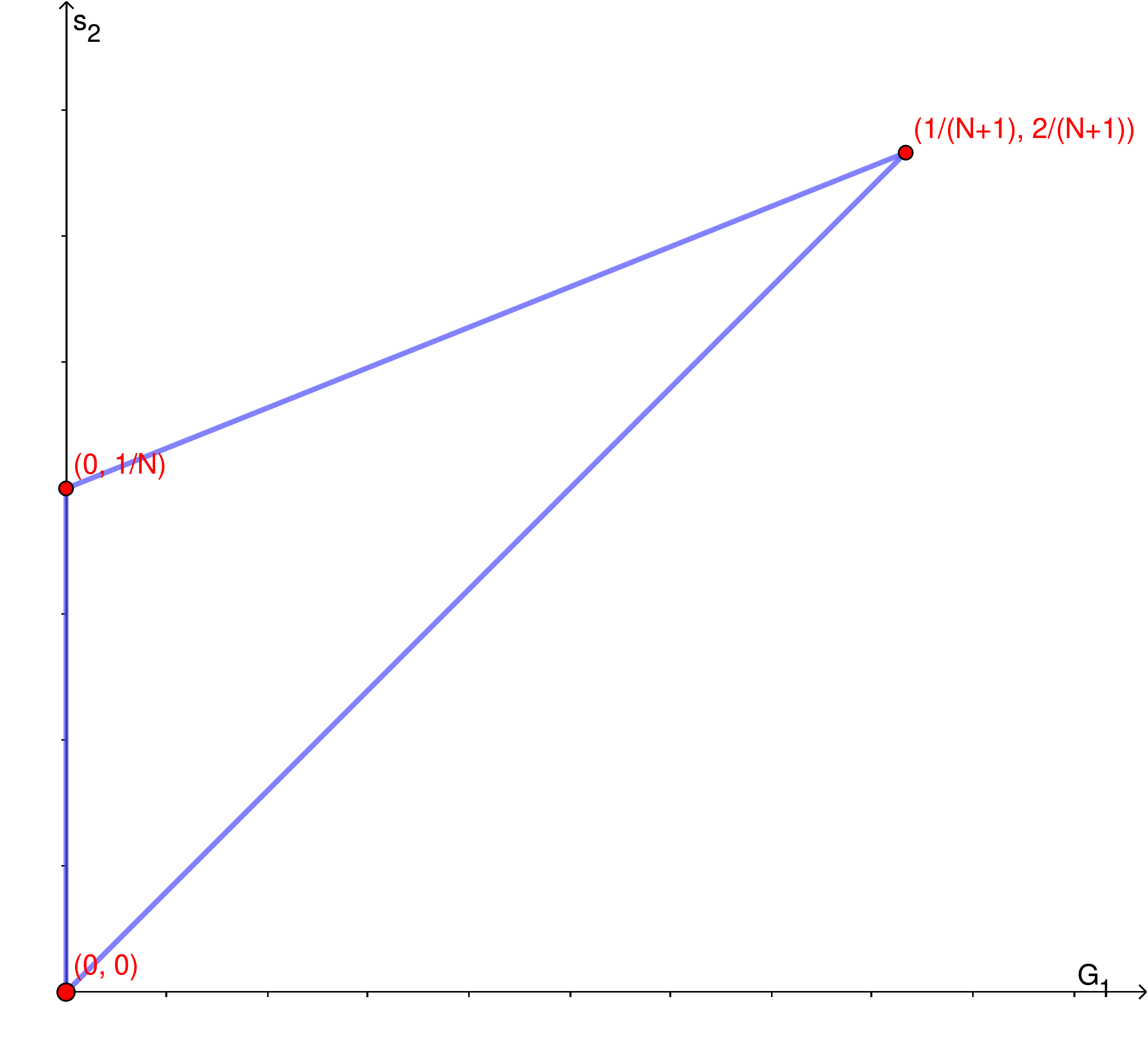}
    \caption{Support of the joint distribution $p(M_1, s_2)$}
    \label{min1_sum2}
\end{figure}

\begin{align} \label{pdf_min2}
    p(s_2 = s) &= \int_{0}^{\frac{s}{2}} (N+1) N^2 (N-1) \left[ 1 + (N - 1)x - N s \right]^{N-2} dx  &&   \text{\small{for $0 \leq s \leq \frac{1}{N}$}} \nonumber \\
    &= \int_{\frac{Ns - 1}{N-1}}^{\frac{s}{2}} (N+1) N^2 (N-1) \left[ 1 + (N - 1)x - N s \right]^{N-2} dx  &&   \text{\small{for $\frac{1}{N} \leq s \leq \frac{2}{N+1}$}} \nonumber \\
    &= \frac{N(N+1)}{(N-1)} \left( \left[ 1 - \left( \frac{N+1}{2} \right)s \right]^{N-1} - \left[ 1 - Ns \right]^{N-1}\right) &&   \text{\small{for $0 \leq s \leq \frac{1}{N}$}} \nonumber \\
    &= \frac{N(N+1)}{(N-1)} \left[ 1 - \left( \frac{N+1}{2} \right)s \right]^{N-1} &&   \text{\small{for $\frac{1}{N} \leq s \leq \frac{2}{N+1}$}}
\end{align}

\subsection{Sum of minima: $k$}

So far we have explicitly derived the joint distributions for the sum of the first two smallest spacings.
We are led to make an hypothesis regarding the general distribution of $s_k$:

\begin{equation} \label{pdf_sumk}
    p(s_k = s | N, 1) = A(k, N) \sum_{i = 1}^{k} a(i, k) \left[ 1 - \left( \frac{N+2-i}{k+1-i} \right) s \right]^{N-1} \cdot H\left( s, 0,  \frac{k + 1 - i}{N + 2 - i} \right) 
\end{equation}

\noindent where the coefficients $A(k,N)$ and $a(i,k)$ are given by:

\begin{equation} \label{coeff_A_sumk}
    A(k, N) = \frac{N (N+1)!}{(N+1-k)^{k-1}(N+1-k)!} 
\end{equation}

\begin{equation} \label{coeff_a_sumk}
    a(i, k) = \frac{(-1)^{i-1}(k+1-i)^{k-2}}{(k-i)!(i-1)!}
\end{equation}
\\
This is indeed the distribution of the sum of the $k$ smallest spacings and a proof by induction is presented in the Appendix. The cumulative density function is:

\begin{align} \label{cdf_sumk}
    p(s_k \leq s | N, 1) &= A(k, N) \sum_{i = 1}^{k} \int_{0}^{\mathrm{min}\left( s, \frac{k + 1 - i}{N + 2 - i} \right) } a(i, k) \left[ 1 - \left( \frac{N+2-i}{k+1-i} \right) x \right]^{N-1} dx \nonumber \\
    &= \frac{A(k, N)}{N} \sum_{i = 1}^{k} \frac{a(i, k)(k+1-i)}{(N+2-i)} \left( 1 - \left[ 1 - \left( \frac{N+2-i}{k+1-i} \right) s \right]^{N} H\left( s,0, \frac{k + 1 - i}{N + 2 - i} \right) \right)
\end{align}

\begin{figure} [h!]
    \centering
    \includegraphics[width=0.95\textwidth]{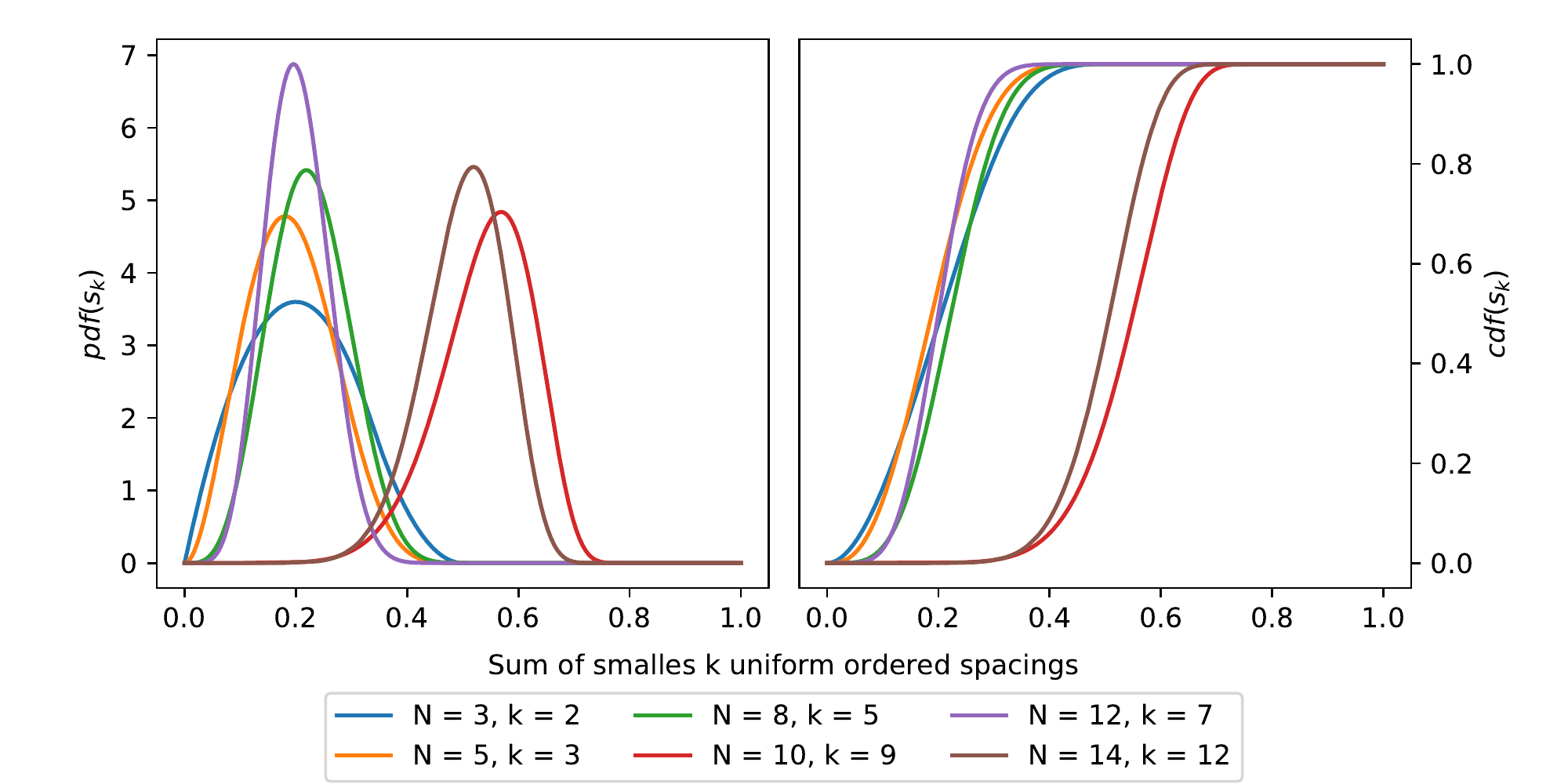}
    \caption{Left: probability distributions for $s_k$, the sum of the $k$ smallest ordered uniform spacings, for different combinations of $N$ and $k$. Right: the cumulative probability distribuion for $s_k$ for the same choices of $N$ and $k$.}
    \label{fig:min1_sum}
\end{figure}

\noindent The distributions that we have derived are plotted in Fig.~\ref{fig:min1_sum} for different choices of $N$ and $k$.

\section{Sum of largest spacings}

\begin{figure} [h!]
    \centering
    \includegraphics[width=0.95\textwidth]{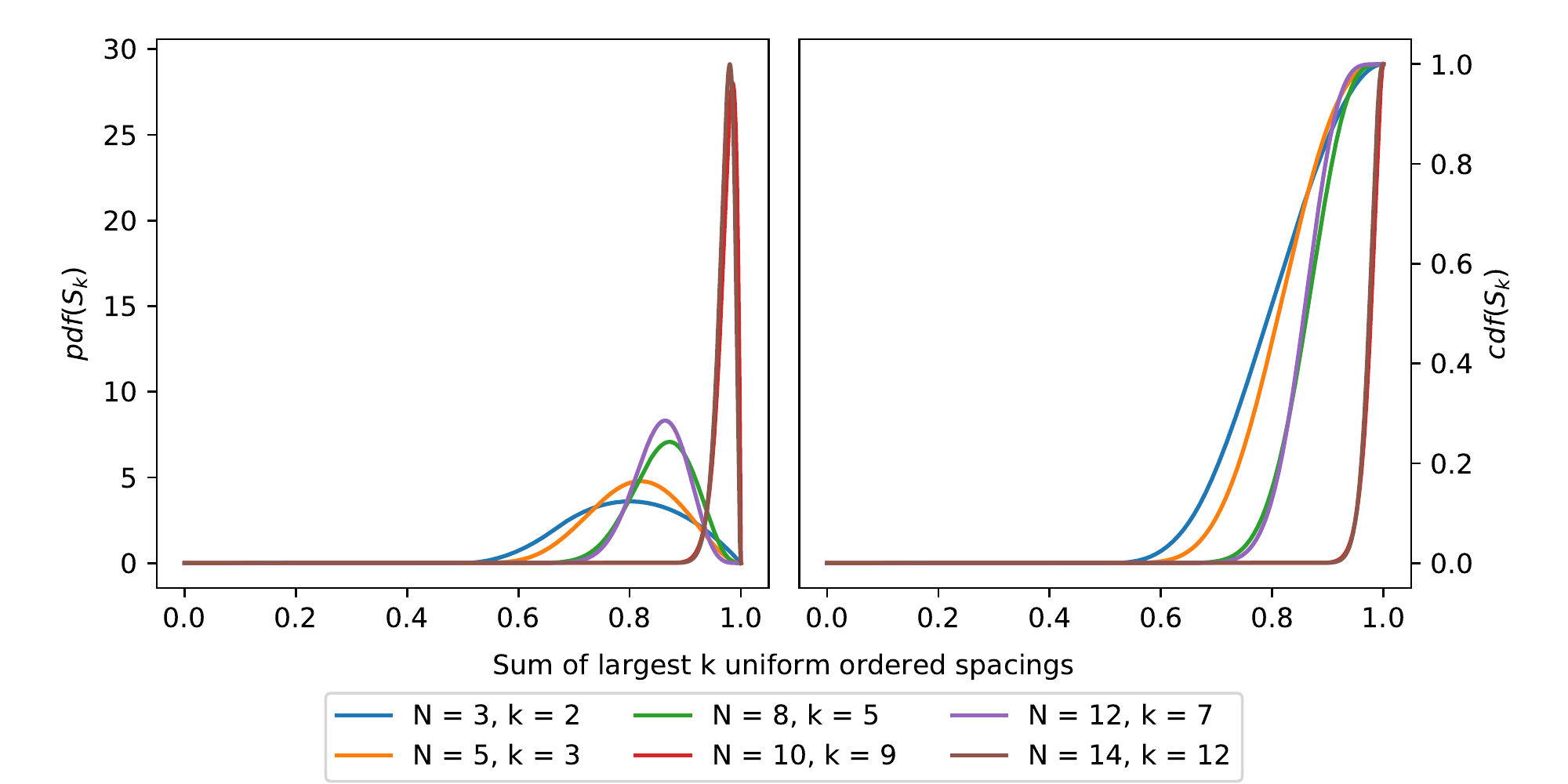}
    \caption{Left: probability distributions for $S_k$, the sum of the $k$ largest ordered uniform spacings, for different combinations of $N$ and $k$. Right: the cumulative probability distribution for $s_k$ for the same choices of $N$ and $k$.}
    \label{fig:max1_sum}
\end{figure}

Since the sum of all the spacings is 1, knowing the sum of the first $k$ smallest spacings allows us to know the value of the sum of the largest $(N+1-k)$ spacings. I.e., we have

\begin{equation}\label{sum_of_maxima}
    S_k = \sum_{i=1}^{k} G_{(N+2-i)} = 1 - \sum_{i=1}^{N+1-k} G_{(i)} = 1 - s_{N+1-k}
\end{equation}
\\
which implies that:

\begin{equation}\label{prob_sum_of_maxima}
    p(S_k = s) = p(s_{N+1-k} = 1-s)
\end{equation}

\begin{equation} \label{pdf_sumk_max}
    p(S_k = s | N, 1) = A(N+1-k, N) \sum_{i = 1}^{N+1-k} a(i, N+1-k) \left[ \frac{s(N+2-i)-k}{N+2-k-i} \right]^{N-1} H\left( s, \frac{k}{N + 2 - i}, 1 \right) 
\end{equation}

\begin{align} \label{cdf_sumk_max}
    p(S_k \geq s | N, 1) &= \frac{A(N+1-k, N)}{N} \sum_{i = 1}^{N+1-k} \frac{a(i, N+1-k)(N+2-k-i)}{N+2-i} \cdot \nonumber \\
    &\quad \cdot \left[ \frac{s(N+2-i)-k}{N+2-k-i} \right]^{N} H\left( s, \frac{k}{N + 2 - i}, 1 \right) 
\end{align}
\\
where the coefficients $A$ and $a$ are the same as Eq.~\ref{coeff_A_sumk}-\ref{coeff_a_sumk}. The distributions that we have derived in plotted in Fig.~\ref{fig:max1_sum} for different choices of $N$ and $k$.

\section{Excluding edges 0 and 1}

So far we have included the boundaries, 0 and 1, in the list of values $X_i$, considering them as de-facto data. This means that in all previous derivations we considered an effective population of $N+2$ values, where $X_{N+1} - X_0 = 1$ and the remaining $N$ values determine the spacings.

\noindent This approach introduces some possibly unwelcome artifacts in the analysis of spacings data. E.g., if we shift all the inner $N$ values slightly towards one of the boundaries? If we were to perform this change, then

\begin{equation} \label{pdf_mink+1_given_mink}
    \{ X_{(1)}, ..., X_{(N)} \} \rightarrow \{ X_{(1)} \pm \epsilon, ..., X_{(N)} \pm \epsilon \} \implies G_1 \pm \epsilon, \quad G_{N+1} \mp \epsilon 
\end{equation}

\noindent Depending on the application, it is possible that one is interested only in the spacings between the inner $N$ values, without considering how close this set is to either boundary.  We therefore would want the distributions that do not include the boundaries. 

In this scenario we will derive the spacing statistic of only the inner spacings $\{ G_2, ..., G_N \}$. In this new scenario, $X_{(1)}$ and $X_{(N)}$ become the new boundaries just as $X_{(0)} = 0$ and $X_{(N+1)} = 1$ were the boundaries in the previous case. 

\begin{equation} \label{pdf_mink+1_given_mink}
    \{ X_{(1)}, ..., X_{(N)} \} \rightarrow \{ X_{(0)}^f, ..., X_{(N-1)}^f \} \rightarrow \{ G_{1}^f, ..., G_{N-1}^f \} \rightarrow \{ G_{(1)}^f, ..., G_{(N-1)}^f \}
\end{equation}
\\
This means that given $N$ values in the no edge scenario we have an effective population of $N-1$ spacings in an interval with length $X_{(N)} - X_{(1)} = X_{(N-1)}^f - X_{(0)}^f = \mu$.\\

\noindent Given $\mu$, we can reuse the same distribution of a quantity we have studied with the presence of the boundaries decreasing the number of spacings from $N+1$ to $N-1$ and by rescaling the support of the distribution to an interval of length $\mu$ by means of the change of variable rule. 
\\
Looking at $\mu$ we notice that it is none other than the spacing between the extremes of the ordered values and its distribution is given by  Eq.~~(\ref{gos_dist}).

Given $N$ values, a quantity of interest $A$ and its distribution with boundaries $p_{w.b.}(A = x | N, \mu)$, in order to derive the distribution of $A$ without boundaries, $p_{n.b.}(A = x | N, 1)$, we have to marginalize over all possible values of $\mu$:

\begin{align} \label{pdf_no_edge}
    p_{n.b.}(A = x | N+1, 1) &= \int_{0}^{1} p(\mu) \cdot p_{w.b.}(A = x | N-1, \mu) d\mu \nonumber \\
    &= \int_{0}^{1} N(N-1)\mu^{N-3}(1-\mu) \cdot p_{w.b.} \left( \left. A = \frac{x}{\mu} \right\vert N-1, 1 \right) d\mu 
\end{align}
\\

\subsection{$k$-th smallest spacing}





The distribution of the $k$-th ordered uniform spacing is given in Eq.~\ref{k-th-spacing}. Using Eq.~\ref{pdf_no_edge} we can get the distribution of the $k$-th ordered uniform spacing without boundaries:

\begin{align} \label{pdf_no_edge_k-th_spacing}
    p_{n.b.}(G_{(k)} = x | N, 1) &=  (-1)^k N (N-1)^2 (N-2) \binom{N-2}{k-1} \cdot \nonumber \\
    &\quad \cdot \sum_{i=1}^k \int_{(N-i)x}^{1} (-1)^i  \binom{k-1}{i-1} (1-\mu) \left[ \mu - (N-i) x \right]^{N-3} d\mu \nonumber \\
    &= (-1)^k N (N-1)^2 \binom{N-2}{k-1} \cdot \nonumber \\
    &\quad \cdot \sum_{i=1}^k \int_{(N-i)x}^{1} (-1)^i  \binom{k-1}{i-1} \left[ \mu - (N-i) x \right]^{N-2} d\mu \nonumber \\
    &= (-1)^k N (N-1) \binom{N-2}{k-1} \cdot \nonumber \\
    &\quad \cdot \sum_{i=1}^k  (-1)^i  \binom{k-1}{i-1} \left[ 1 - (N-i) x \right]^{N-1}  H \left(x, 0, \frac{1}{N-i} \right) 
\end{align}

\begin{align} \label{cdf_no_edge_k-th_spacing}
    p_{n.b.}(G_{(k)} \leq x | N, 1) &= (-1)^k (N-1) \binom{N-2}{k-1} \cdot \nonumber \\
    &\quad \cdot \sum_{i=1}^k  \frac{(-1)^i}{(N-i)}  \binom{k-1}{i-1} \left( 1 - \left[ 1 - (N-i) x \right]^{N}  H \left(x, 0, \frac{1}{N-i} \right) \right)
\end{align}

\begin{figure} [h!]
    \centering
    \includegraphics[width=0.95\textwidth]{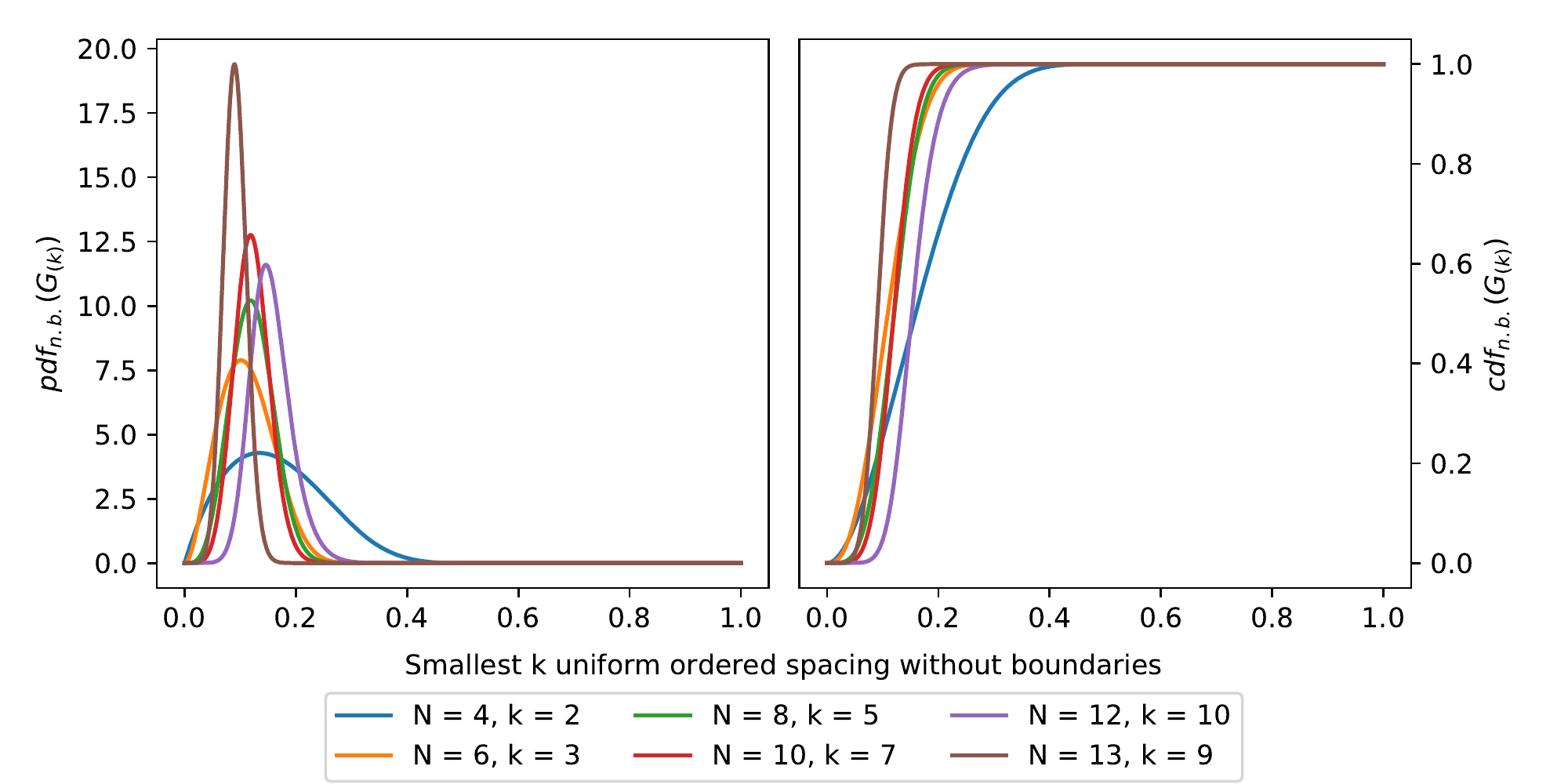}
    \caption{Left: probability distributions for $G(k)$, the $k^{\rm th}$ smallest spacing, for different combinations of $N$ and $k$, for the case where the boundaries are not included in the definition of the spacings. Right: the cumulative probability distribuion for $G(k)$ for the same choices of $N$ and $k$.}
    \label{min2_nosum}
\end{figure}

Examples of the resulting distributions are shown in Fig.~\ref{min2_nosum}.  



\subsection{Sum of $k$ smallest spacings}

For the sum of the first $k$ uniform ordered spacings we get:

\begin{align} \label{pdf_no_edge_sum_min_k_spacings}
    p_{n.b.}(s_k = s | N, 1) &=  N (N-1) A(k, N-2) \sum_{i = 1}^{k} \int_{\frac{(N-i)s}{k+1-i}}^{1} a(i, k) (1 - \mu) \left[ \mu - \left( \frac{N-i}{k+1-i} \right) s \right]^{N-3} d\mu \nonumber \\
    &= \frac{N}{(N-2)} A(k, N-2) \sum_{i = 1}^{k}  a(i, k) \left[ 1 - \left( \frac{N-i}{k+1-i} \right) s \right]^{N-1} H\left( s, 0, \frac{k+1-i}{N-i} \right)  
\end{align}

\begin{align} \label{cdf_no_edge_sum_min_k_spacings}
    p_{n.b.}(s_k \leq s | N, 1) &= \frac{A(k, N-2)}{(N-2)}  \sum_{i = 1}^{k}  \frac{a(i, k)(k+1-i)}{(N-i)} \left( 1 -\left[ 1 - \left( \frac{N-i}{k+1-i} \right) s \right]^{N} H\left( s, 0, \frac{k+1-i}{N-i} \right)  \right)
\end{align}

\begin{figure} [h!]
    \centering
    \includegraphics[width=0.95\textwidth]{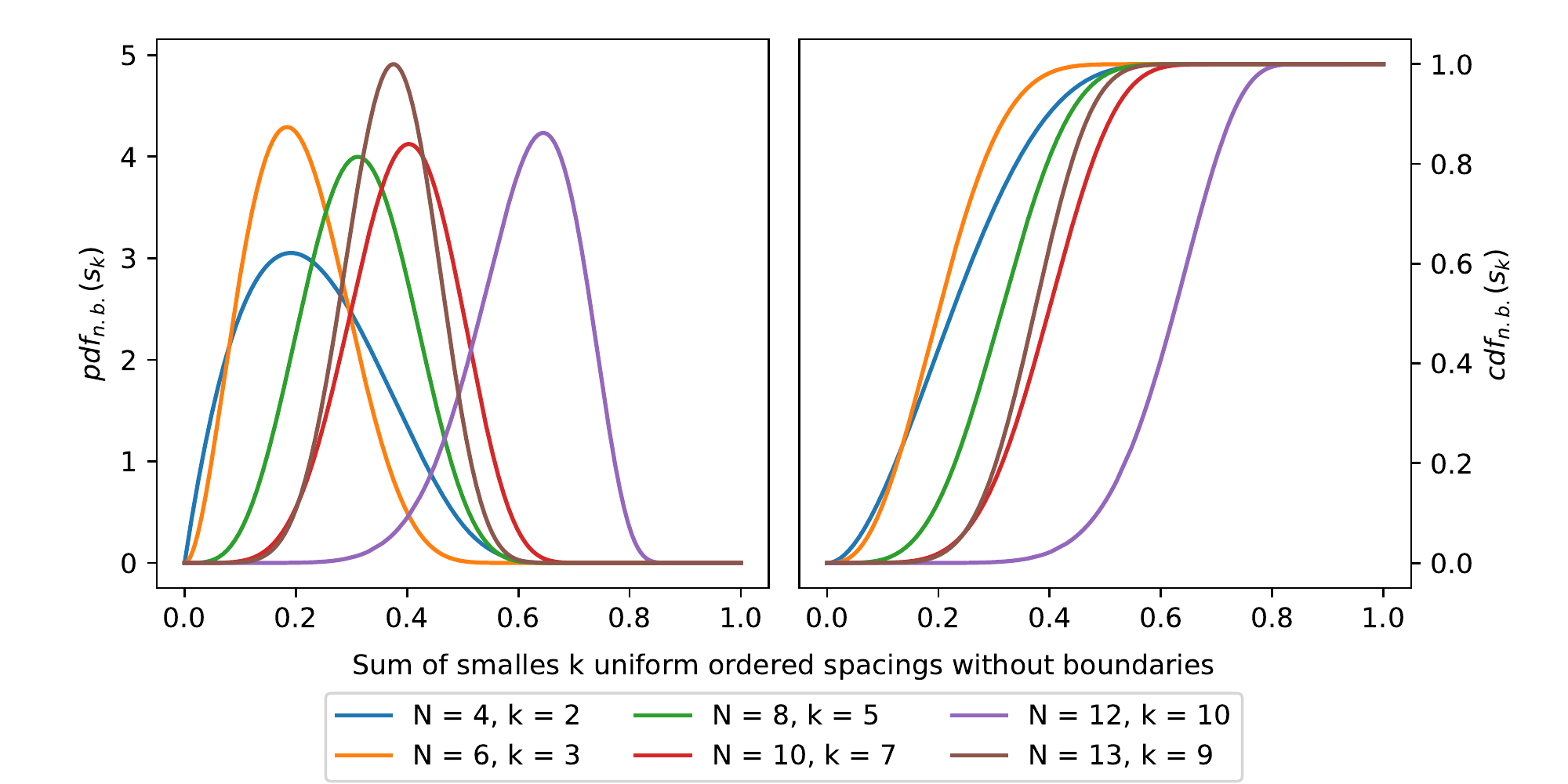}
    \caption{Left: probability distributions for $s_k$, the sum of the $k$ smallest ordered uniform spacings, for different combinations of $N$ and $k$ for the case where the boundaries are not included in the definition of the spacings. Right: the cumulative probability distribuion for $s_k$ for the same choices of $N$ and $k$.}
    \label{min2_sum}
\end{figure}

Our expressions are displayed in Fig.~\ref{min2_sum}.

\subsection{Sum of $k$ largest spacings}

\begin{figure} [h!]
    \centering
    \includegraphics[width=0.95\textwidth]{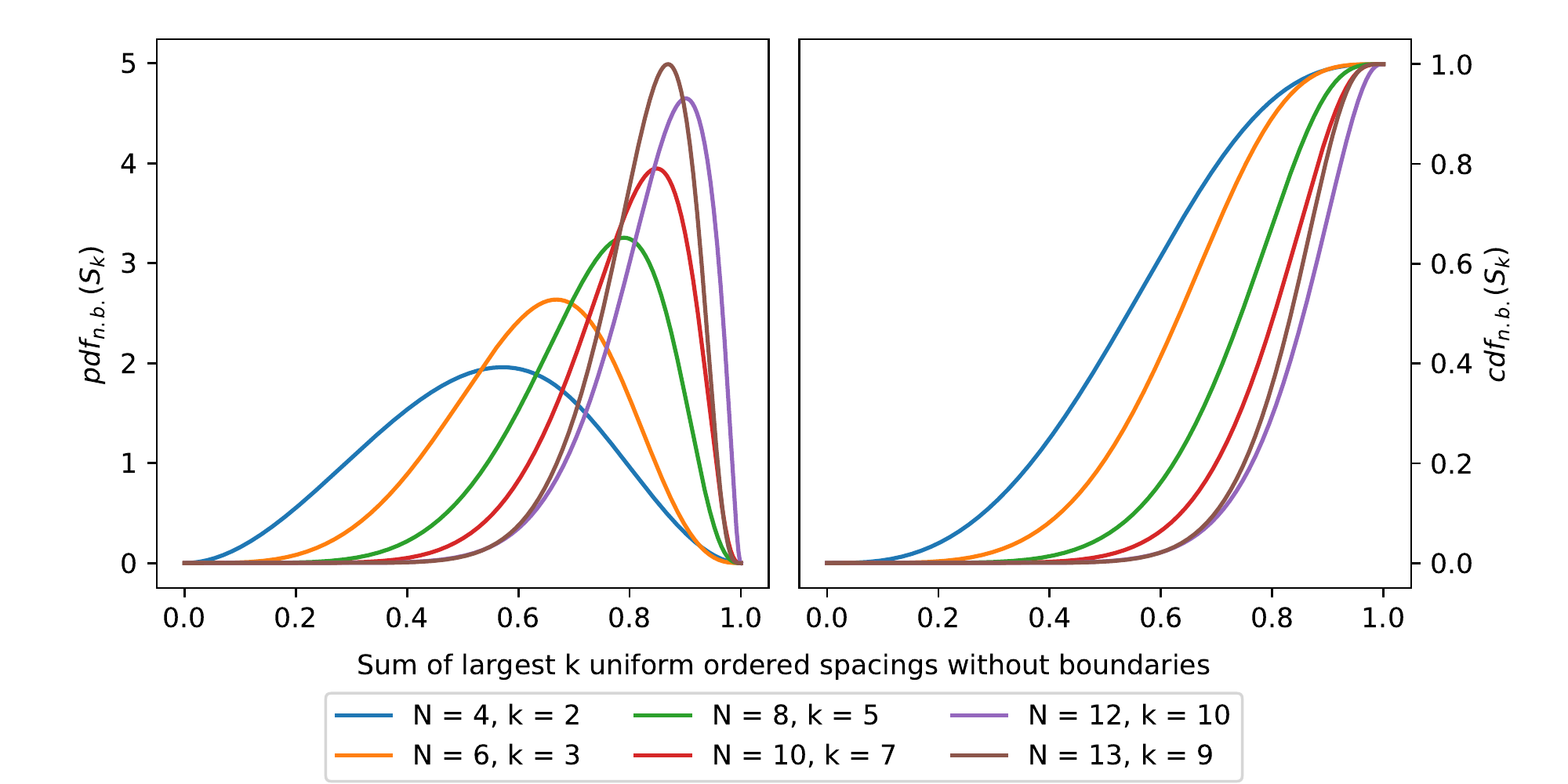}
    \caption{Left: probability distributions for $S_k$, the sum of the $k$ largest ordered uniform spacings, for different coombinations of $N$ and $k$ for the case where the boundaries are not included in the definition of the spacings. Right: the cumulative probability distribuion for $s_k$ for the same choices of $N$ and $k$.}
    \label{max2_sum}
\end{figure}

For the sum of the last $k$ uniform ordered spacings we get:

\begin{align} \label{pdf_no_edge_sum_max_k_spacings}
    p_{n.b.}(S_k = s | N, 1) &=  N (N-1) A(N-1-k, N-2) \cdot \nonumber \\
    &\quad \cdot \sum_{i = 1}^{N-1-k} \int_{s}^{\mathrm{min} \left( 1, \frac{s(N-1)}{k} \right)}  a(i, N-1-k) (1 - \mu) \left[ \frac{s(N-i) - \mu k}{N-k-i} \right]^{N-3} d\mu \nonumber \\
    &= \frac{N A(N-1-k, N-2)}{k^2(N-2)}  \sum_{i = 1}^{N-1-k} \frac{a(i, N-1-k)}{(N-k-i)^{N-3}} \cdot \left( \right. \nonumber \\
    &\quad  s^{N-2} \left[ k(N-1) - s(N-i) -ks(N-2) \right] [N-i-k]^{N-2}  +  \nonumber \\
    &\quad \left. + [s(N-i) - k]^{N-1}H \left( s, \frac{k}{N-i}, 1 \right) \right)
\end{align}

\begin{align} \label{cdf_no_edge_sum_max_k_spacings}
    p_{n.b.}(S_k \leq s | N, 1) &= \frac{N A(N-1-k, N-2)}{k^2(N-2)}  \sum_{i = 1}^{N-1-k} \frac{a(i, N-1-k)}{(N-k-i)^{N-3}} \cdot \left( \right. \nonumber \\
    &\quad  s^{N-1} \left[k - \frac{x \left[ N(k + 1) - i - 2k \right]}{N} \right] [N-i-k]^{N-2}  +  \nonumber \\
    &\quad \left. + \frac{1}{N(N - i)} [s(N-i) - k]^{N}H \left( s, \frac{k}{N-i}, 1 \right) \right)
\end{align}
\\
Examples of these are displayed Fig.~\ref{max2_sum}

\section{Discussion}

We have derived the probability distributions of the sums of either the smallest or the largest $k$ ordered uniform spacings. We did this both for the case where the boundaries are included in the analysis and when only the observed values are used.

These quantities can be very useful when analysing a sequence of results in a particle physics context.  They can either give an indication for the presence of a unexpected source via a clustering of values, or they can be used to set an upper limit on the normalization of a spectrum.   We note that in many experimental scenarios, the number of observed values is a random variable.  The distributions we have derived can then be used convoluted with the expected distribution for the random number of observed values.

We are currently developing novel test statistics for this purpose for use in particle physics data analysis.

\printbibliography

\newpage

\section*{Appendix}

In the previous sections we have shown that Eq.\ref{pdf_sumk} is valid for $s_1$ and $s_2$. We now show its validity going from $k-1$ to $k$ via induction. We start from the joint distribution of $s_k$ and $G_{(1)}$:

\begin{align} \label{pdf_sumk_bayes}
    p(s_k &= s, G_{(1)} = x | N, 1) = p(G_{(1)} = x| N, 1)p(s^*_{k-1} = s - kx | N-1, 1 - (N+1)x) \nonumber \\
    &= p(G_{(1)} = x | N, 1) \left( \frac{1}{1 - (N+1)x} \right) p\left(s_{k-1}^* = \frac{s - kx}{1 - (N+1)x} | N-1, 1 \right) \nonumber \\
    &= N(N+1) \left[ 1 - (N+1)x \right]^{N-2} A(k-1, N-1) \cdot \left( \sum_{i = 1}^{k-1} a(i, k-1) \cdot  \right. \nonumber \\
    &\quad \left. \cdot \left[ 1 - \left( \frac{N+1-i}{k-i} \right) \frac{s - kx}{1 - (N+1)x}  \right]^{N-2} \cdot H\left( \frac{s - kx}{1 - (N+1)x}, 0, \frac{k - i}{N + 1 - i} \right) \right) \nonumber \\ 
    &= N(N+1) A(k-1, N-1) \cdot \left( \sum_{i = 1}^{k-1} a(i, k-1) \cdot  \right. \nonumber \\
    &\quad \left. \cdot \left[ 1 + x \cdot \frac{i(N+1-k)}{k-i} - s \cdot \frac{(N+1-i)}{k-i} \right]^{N-2} \cdot H\left( \frac{s - kx}{1 - (N+1)x}, 0, \frac{k - i}{N + 1 - i} \right) \right)
\end{align}

\noindent Marginalizing over $M_1$ we have:

\begin{align} \label{pdf_sumk_int}
    p(s_k &= s| N, 1) = \int_{0}^{\frac{s}{k}} p(s_k = s, M_1 = x | N, 1) \nonumber \\
    &= N(N+1) A(k-1, N-1) \cdot \nonumber \\
    &\quad \cdot \sum_{i = 1}^{k-1} \int_{\textrm{max}\left( 0, \frac{s(N+1-i) - k + i}{i(N+1-k)} \right)}^{\frac{s}{k}} a(i, k-1) \left[ 1 + x \cdot \frac{i(N+1-k)}{k-i} - s \cdot \frac{(N+1-i)}{k-i} \right]^{N-2}  \nonumber \\
    &= \frac{N(N+1)}{(N-1)(N+1-k)} A(k-1, N-1) \cdot \nonumber \\
    &\quad \cdot \sum_{i = 1}^{k-1}  a(i, k-1) \cdot \frac{(k-i)}{i} \cdot \left. \left( \left[ 1 + x \cdot \frac{i(N+1-k)}{k-i} - s \cdot \frac{(N+1-i)}{k-i} \right]^{N-1}  \right\vert_{\textrm{max}\left( 0, \frac{s(N+1-i) - k + i}{i(N+1-k)} \right)}^{\frac{s}{k}} \right) \nonumber \\
    &= \frac{N(N+1)}{(N-1)(N+1-k)} A(k-1, N-1) \cdot \sum_{i = 1}^{k-1}  a(i, k-1) \cdot \frac{(k-i)}{i} \cdot  \left( \left[ 1 - s \cdot \frac{(N+1)}{k} \right]^{N-1} \right. - \nonumber \\
    &\quad - \left. \left[ 1 - s \cdot \frac{(N+1-i)}{k-i} \right]^{N-1} H\left( s, 0, \frac{k-i}{N+1-i} \right)  \right) \nonumber \\
    &= \frac{N(N+1)}{(N-1)(N+1-k)} A(k-1, N-1) \cdot \left( \left( \sum_{i = 1}^{k-1}  a(i, k-1) \cdot \frac{(k-i)}{i} \right) \cdot   \left[ 1 - s \cdot \frac{(N+1)}{k} \right]^{N-1} \right. - \nonumber \\
    &\quad - \left. \sum_{i = 2}^{k}  a(i-1, k-1) \cdot \frac{(k+1-i)}{i-1} \cdot \left[ 1 - s \cdot \frac{(N+2-i)}{k+1-i} \right]^{N-1} \cdot H\left( s, 0, \frac{k+1-i}{N+2-i} \right)  \right)
\end{align}
\\
Looking back at Eq.~\ref{coeff_A_sumk} notice that:

\begin{align} \label{update_A_coeff}
    \frac{N(N+1)}{(N-1)(N+1-k)} A(k-1, N-1) &= \frac{N(N+1)}{(N-1)(N+1-k)} \cdot \frac{(N-1) N!}{(N-k-1)^{k-2}(N-k-1)!} \nonumber \\
    &= \frac{N (N+1)!}{(N+1-k)^{k-1}(N+1-k)!} \nonumber \\
    &= A(k, N)
\end{align}

\begin{align} \label{update_a_coeff}
    -a(i-1, k-1) \cdot \frac{(k+1-i)}{i-1} &= - \frac{(-1)^{i-2}(k+1-i)^{k-3}}{(k-i)!(i-2)!} \cdot \frac{(k+1-i)}{i-1} \nonumber \\
    &= \frac{(-1)^{i-1}(k+1-i)^{k-2}}{(k-i)!(i-1)!} \nonumber \\
    &= a(i, k)  &&   \text{\small{for $2 \leq i \leq k$}}
\end{align}

\noindent The result of Eq.~\ref{update_a_coeff} implies a recursion formula for the coefficients $a(i, k) = f[a(i-1, k-1)]$. Making use of this recursion we can relate any $a(i, k)$ to $a(1, k+1-i)$:

\begin{equation} \label{update_a_i_to_a_1}
    a(i, k) = \frac{(-1)^{i-1}(k+1-i)^{i-1}a(1, k+1-i)}{(i-1)!}
\end{equation}
\\
Finally we have that:

\begin{align} \label{update_a_1_coeff_recursion}
    \sum_{i = 1}^{k-1}  a(i, k-1) \cdot \frac{(k-i)}{i} &= \sum_{i = 1}^{k-1}  \frac{(-1)^{i-1}(k-i)^{i-1}a(1, k-i)}{(i-1)!} \cdot \frac{(k-i)}{i} \nonumber \\
    &= - \sum_{i = 1}^{k-1}  \frac{(-1)^{i}(k-i)^{i}a(1, k-i)}{i!} 
\end{align}
\\
In order for Eq.~\ref{pdf_sumk_int} to satisfy our hypothesis, we need that:

\begin{equation} \label{sum_a_1_coeff}
    \sum_{i = 1}^{k-1}  a(i, k-1) \cdot \frac{(k-i)}{i} = a(1, k)
\end{equation}
\\
Putting together Eq.~\ref{update_a_1_coeff_recursion} and Eq.~\ref{sum_a_1_coeff} we find a recursion rule for the coefficients of the coefficients $a(1, k)$. Using this recursion we get:

\begin{align} \label{update_a_1_coeff}
    - \sum_{i = 1}^{k-1} &\quad \frac{(-1)^{i}(k-i)^{i} a(1, k-i)}{i!} = - \sum_{i = 1}^{k-2} \frac{ik}{(i+1)} \cdot \frac{(-1)^{i}(k-1-i)^{i}a(1, k-1-i)}{i!}  \nonumber \\
    &= - \sum_{i = 1}^{k-m} \frac{i \cdot k^{m-1} (-1)^{i}(k-m+1-i)^{i}a(1, k-m+1-i)}{(m-1)!(i+m-1) i!} \text{\small{    for $1 \leq m \leq k-1$}} \nonumber \\
    &= \frac{k^{k-2}}{(k-1)!} \cdot a(1,1)  =  a(1, k)  
\end{align}
\\
where we have used Eq.~\ref{sum_a_1_coeff} to express the first factor in each of the sums, allowing us to reduce the limits of the sum by means of this recursion. The result we obtain proves the consistency of the recursion relation which finally proves the consistency of our hypothesis.

\end{document}